\documentclass[11pt,leqno]{article}

\usepackage{amsmath} 
\usepackage{amsfonts} 
\usepackage{amssymb} 
\newtheorem{theo}{Theorem} 
\usepackage[latin1]{inputenc} 
\usepackage[hyperfootnotes=false,colorlinks]{hyperref}

\begin{document}

\title{Some convolution inversion questions} 
\date{14 mai 2019} 
\author{Michel {\sc Valadier}} 
\maketitle

\begin{abstract} 
Blurring of a photographic image by a wrong focus 
can be modeled by convolution. 
Is inversion a possible answer? 
This paper adds complements to a foregoing paper \cite{Va} 
discussing convolution-inversion of some measures. 
\end{abstract} 


\section{Introduction}\label{Int} 
\noindent 
The inverse operation of the ``blurring'' effect induced 
by convolution (for example wrong focus, cf.~\cite[Introduction]{Va}) 
seems not too badly achieved 
by a \emph{high-pass filter} 
i.e.\ by convolution with a kernel as \emph{the mexican 
hat}\footnote{\ In many documents, too often, 
this is shown with a $3 \times 3$ or $5 \times 5$ matrix of pixels!} 
(or with more oscillations; cf.\ the \emph{cardinal sine function} 
which infinitely oscillates). 
My naive explanation: 
using linear mixing\footnote{\ With coefficients whose algebraic 
sum is $1$ in order to keep the same mean value.} of values 
of the blurred signal seems natural; and this using only values at 
neighbouring points; 
and with some negative coefficients, otherwise this would increase 
the blurring\footnote{\ This is rather qualitative: 
I will not give values of amplitude 
and wavelength for any mexican hat.}. 
As for values at distant points, only academical examples 
\cite[equation~(10) and Section~5 ``Case of almost perfect grey'']{Va} 
could use them. 

This is a small side of ``Image Processing''. 
For a large point of view see \cite{J} (book for engineers?), 
for a quasi-philosophical view see \cite{M2} (maybe 
\cite{M1} is a shorter version) 
and there are abundance of papers using wavelets (Ingrid {\sc Daubechies} 
being surely the most famous author; her most quoted work: \cite{Da}). 
In all the literature these problems are said \emph{ill-posed}. 
For a now old book cf.~\cite[Ch.IV pp.91--115]{TA}. 

We don't bring any general solution but continue the observations 
of \cite{Va} about prospective inverses. 
We recall some definitions in Section~\ref{inver}, prove an inversion 
result in Section~\ref{easy} where the value $1/2$ plays a central role, 
and give a strange inverse in Section~\ref{calc} 
(completing thus \cite[Th.3.6]{Va}). 
Finally we give some naive observations about 
the gaussian kernel in Section~\ref{gauss}. 
\medskip

\section{About inverses}\label{inver} 
\noindent 
Let $T$ be a distribution on $\mathbb R^d$ 
(about distributions an historical reference is \cite{S}). 
Its \emph{convolution product} with any other distribution 
makes sense if $T$ has compact support but also in many 
different cases. 
The unit element is $\delta_0$ (the measure with mass $1$ 
at the point $0$). 
The distributions $T$ and $D$ constitute a 
\emph{couple of zero divisors} if $T \neq 0$, $D \neq 0$ and 
\begin{equation}\label{div} 
T \ast D = 0 \,. 
\end{equation} 
And $V$ (this for $T$ necessarily $\neq 0$) is an \emph{inverse} 
of $T$ if 
\begin{equation}\label{inv} 
T \ast V = \delta_0 \,, 
\end{equation} 
in which case also $V \neq 0$. 
The set possibly empty of all inverses of $T$ could 
be denoted by $T^{\ast(-1)}$. 
With this notation \eqref{inv} writes 
\[ 
V \in T^{\ast(-1)} \,. 
\] 
If $V_1$ and $V_2$ are inverses of $T$, for any $\lambda \in \mathbb R$, 
$\lambda \, V_1 + (1 - \lambda) V_2$ is also an inverse of $T$.
If \eqref{inv} and \eqref{div} hold, $V + D$ is still an inverse of $T$. 
If $V_1$ and $V_2$ are two different inverses of $T$, 
then $T$ and $V_2 - V_1$ 
constitute a couple of zero divisors. 
We gave examples in \cite[Section~2]{Va}. 

One may ask if there exist distributions which are not 
zero divisors. 
Easy answer: $\delta_x$ for any $x \in \mathbb R^d$. 
When associativity holds any distribution admitting an inverse 
is not a zero divisor. 

One could ask also: does there exist distributions with 
not any inverse? 
In my opinion surely: see Section~\ref{gauss}. 

Surely spaces of distributions (or of general measures) are too large. 
Images have compact supports! 
They have density and even more should be considered only as pixels. 
And the Fourier transform with the ``temp\'er\'ees distributions'' 
(French terminology) is too magical: cf.~the Dirac comb. 
Precise frameworks would be essential. 
\medskip

\section{An easy inversion}\label{easy} 
\noindent 
For $x \in \mathbb R$, $|x| <1$: 
\[ 
(1 + x)^{-1} = 1 - x + x^2 - x^3 + \dots \,, 
\] 
and this series seems the idea of Van Cittert\footnote{\ See 
in German \cite{VC2}. 
This Wikipedia page quotes B.~J\"ahne whose book is now 
in its sixth edition; I don't see it. 
See specially, in the fifth edition \cite{J} 17.8.4 pp.\,478--480. 
Reference \cite{VC1} is given by \cite{B} and in latest 
versions of Wikipedia.}. 

Let $d$ be an integer $\geq 1$, 
$\mathcal M^b(\mathbb R^d)$ or shortly $\mathcal M^b$ the set 
of bounded measures on $\mathbb R^d$. 
The \emph{total variation norm} of $\mu \in \mathcal M^b$ is 
\[ 
\|\mu\| = \mu^+(\mathbb R^d) + \mu^-(\mathbb R^d) \,. 
\] 
Denoting $\mathcal C_k(\mathbb R^d)$ 
the space or continuous functions with compact support on $\mathbb R^d$, 
and $\|\varphi\|$ the uniform norm of 
$\varphi \in \mathcal C_k(\mathbb R^d)$, there holds 
\[ 
\|\mu\| = \sup\Bigl\{\int \varphi \, d\mu \,;\, 
\varphi \in \mathcal C_k(\mathbb R^d), \, \|\varphi\| \leq 1\Bigr\} \,. 
\] 
Equipped with the convolution product, $\mathcal M^b$ is a 
Banach algebra. 
For $n \geq 1$, we note $\mu^{n\ast}$ the $n$-power 
convolution of $\mu$ i.e. 
\[ 
\mu^{n\ast} = \underbrace{\mu \ast \dots \ast \mu}_{n} \,. 
\]

\begin{theo}\label{T1} 
Let $\mu \in \mathcal M^b(\mathbb R^d)$ satisfying $\|\mu\| < 1$. 
Then the series 
\begin{equation}\label{serie} 
\delta_0 + \sum_{k = 1}^{\infty} (-1)^k \, \mu^{k\ast} 
\end{equation} 
converges in $\mathcal M^b$ to a measure $\nu$ 
which is a convolution inverse of $\delta_0 + \mu$ in the space of 
bounded measures on $\mathbb R^d$. 
\end{theo} 

\noindent 
{\sc Remark.} 
If $\mu \geq 0$ (and still $\|\mu\| < 1$) and one considers 
\[ 
\frac{1}{1 + \|\mu\|} (\delta_0 + \mu) 
\] 
we have, as in Duval \cite{Du}, a probability measure where $\delta_0$ 
weight is $> \frac{1}{2}$ (note that Duval, besides this 
hypothesis, uses wavelets!). 
\medskip 

\noindent 
{\sc Proof.} 
For two bounded measures $\lambda_1$ and $\lambda_2$, there holds 
\[ 
\|\lambda_1 \ast \lambda_2 \| \leq \|\lambda_1\| \, \|\lambda_2\| 
\] 
because, for $\varphi \in \mathcal C_k(\mathbb R^d)$ with norm $\leq 1$, 
one has 
\begin{align*} 
\int_{\mathbb R^d} \varphi \, d(\lambda_1 \ast \lambda_2) 
&= \int_{\mathbb R^d} 
\Bigl[ \int_{\mathbb R^d} \varphi(x + y) \, d\lambda_1(x) \Bigr] 
 \, d\lambda_2(y) \\ 
& \leq \int_{\mathbb R^d} \|\lambda_1\| \, d\lambda_2(y) 
\leq \|\lambda_1\| \,\|\lambda_2\| \,. 
\end{align*} 
Hence $\|\mu^{k\ast}\| \leq \|\mu\|^k$. 
Let us set 
\[ 
\nu_n = \delta_0 + \sum_{k = 1}^n (-1)^k \, \mu^{k\ast} \, .
\] 
This series does converge (it is absolutely convergent) and, as 
\[ 
(\delta_0 + \mu) \ast \nu_n = \delta_0 + (-1)^n\,\mu^{(n+1)\ast} 
\] 
there holds 
\[ 
(\delta_0 + \mu) \ast \lim_n \nu_n = \delta_0 \,. \quad\Box 
\] 
\medskip 

In \cite[Theorem~3.6]{Va} an inverse of 
$\frac{1-a}{2} \, \delta_{-1} + a \, \delta_0 
+ \frac{1-a}{2} \, \delta_1$, 
with the parameter $a$ belonging to $\left]\frac{1}{2},1\right[$, 
is given. 
Surely Theorem~\ref{T1} could apply because 
\[ 
\frac{1-a}{2} \, \delta_{-1} + a \, \delta_0 
+ \frac{1-a}{2} \, \delta_1 = 
a \, \Bigl[\delta_0 + \frac{1-a}{2a} \, (\delta_{-1} + \delta_1)\Bigr] 
\] 
and, thanks to $a \in \left]\frac{1}{2},1\right[$, 
one has 
$\bigl\|\frac{1-a}{2a} \, (\delta_{-1} + \delta_1)\bigr\| < 1$. 
Thus setting $\mu = \frac{1-a}{2a} \, (\delta_{-1} + \delta_1)$, 
hypotheses of Theorem~\ref{T1} are verified.

\section{A strange inverse}\label{calc} 
The measure 
\begin{equation}\label{onehalf} 
\mu := 
\frac{1}{4} \, \delta_{-1} + \frac{1}{2} \, \delta_0 
+ \frac{1}{4} \, \delta_1 = 
\frac{1}{2} \, \Bigl[\delta_0 + \frac{1}{2} \, (\delta_{-1} 
+ \delta_1) \Bigr] 
\end{equation} 
is the limit in Theorem~3.6 of \cite{Va} when $a$ tends to $1/2$: 
there the coefficients in the expression of the inverse 
explode to infinity. 
And we leaved open the question of existence of an inverse.

But another way is possible. 
``Lateral'' inversion in $\mathcal D'_+$ (or in $\mathcal D'_-$) 
is possible. 
Indeed 
\[ 
\frac{1}{4} \, \delta_{-1} + \frac{1}{2} \, \delta_0 
+ \frac{1}{4} \, \delta_1 = \frac{1}{4}
(\delta_{-1} + \delta_0) \ast (\delta_0 + \delta_1) 
\] 
and (cf.~\cite[Lemma~3.1]{Va}) 
$\delta_0 + \delta_1$ admits as inverse 
\begin{equation}\label{first} 
\delta_0 - \delta_1 + \delta_2 - \delta_3 + \dots 
\end{equation} 
just as $\delta_{-1} + \delta_0$ admits as inverse 
\begin{equation}\label{second} 
\delta_1 - \delta_2 + \delta_3 - \delta_4 + \dots 
\end{equation} 
Multiplying term by term \eqref{first} and \eqref{second} gives 
\begin{equation}\label{dr} 
\Bigl[(\delta_{-1} + \delta_0) \ast (\delta_0 + \delta_1)\Bigr]^{*(-1)} 
\ni \delta_1 - 2 \, \delta_2 + 3 \, \delta_3 - 4 \, \delta_4 + \dots 
\end{equation} 
(and multiplying this by $4$ would give an inverse 
of $\mu$ in $\mathcal D'_+$). 
But this can also be done on left, in $\mathcal D'_-$: 
an inverse of $\delta_0 + \delta_1$ is 
\begin{equation}\label{first*} 
\delta_{-1} - \delta_{-2} + \delta_{-3} + \dots 
\end{equation} 
while $\delta_{-1} + \delta_0$ admits as inverse 
\begin{equation}\label{second*} 
\delta_0 - \delta_{-1} + \delta_{-2} - \delta_{-3} + \dots 
\end{equation} 
There multiplying \eqref{first*} and \eqref{second*} gives 
\begin{equation}\label{g} 
\Bigl[(\delta_{-1} + \delta_0) \ast (\delta_0 + \delta_1)\Bigr]^{*(-1)} 
\ni \delta_{-1} - 2 \, \delta_{-2} + 3 \, \delta_{-3} - 4 \, \delta_{-4} 
+ \dots 
\end{equation} 
and taking the half sum of the right members of \eqref{dr} and 
\eqref{g} one gets the following symmetric inverse of 
$(\delta_{-1} + \delta_0) \ast (\delta_0 + \delta_1)$ 
\[ 
\ldots + \frac{3}{2} \, \delta_{-3} - \frac{2}{2} \, \delta_{-2} 
+ \frac{1}{2} \, \delta_{-1} + \frac{1}{2} \, \delta_{1} 
- \frac{2}{2} \, \delta_{2} + \frac{3}{2} \, \delta_{3} + \dots 
\] 
This can be checked directly. 
Note that the coefficients oscillate, take negative values 
and tend to infinity. 
Multiplying by $4$ we get the following

\begin{theo}\label{T2} 
The measure $\mu$ defined by \eqref{onehalf} admits as inverse 
in $\mathcal D'(\mathbb R)$ (or also in the space of all measures 
on $\mathbb Z$) 
\[ 
\nu := 
\ldots + 6 \, \delta_{-3} - 4 \, \delta_{-2} 
+ 2 \, \delta_{-1} + 2 \, \delta_{1} 
- 4 \, \delta_{2} + 6 \, \delta_{3} + \dots 
\] 
where for $n \in \mathbb Z$ the coefficient of $\delta_n$ is 
$2\,|n|\,(-1)^{|n|+1}$. 
Thus for any image\footnote{\ We denote as Bourbaki by 
$\mathbb R^{(\mathbb Z)}$ the set of all sequences on $\mathbb Z$ 
with compact support.} $f \in \mathbb R^{(\mathbb Z)}$ 
with compact support there holds $(f \ast \mu) \ast \nu = f$. 
\end{theo} 

As pessimistic (with respect to treatment of badly focused images) 
observations note that in Theorem~3.3 of \cite{Va} 
we got for $1/2 \, (\delta_0 + \delta_1)$ the inverse 
\[ 
H = \ldots - \delta_{-4} + \delta_{-3} - \delta_{-2} + \delta_{-1} 
+ \delta_0 - \delta_1 + \delta_2 - \delta_3 + \ldots 
\] 
whose coefficients do not tend to $0$, and that in Theorem~2 
above the coefficients tend to $+\infty$.

\section{The gaussian kernel case}\label{gauss} 
Let $f$ be a real valued function on $\mathbb R^d$ 
belonging to $L^1$ or even to $L^1 \cap L^2$. 
We are interested by $f$ (for example the luminous intensity 
of a monochrom image; this could be written with $d=2$) 
but observe $g$ given by 
\begin{equation}\label{blur} 
g = f \ast h 
\quad \text{where} \quad 
h(x) = \frac{1}{\bigl(\sqrt{2 \pi}\,\bigr)^d} \, 
e^{-\|x\|^2/2} \,. 
\end{equation} 
Note that $h$ is the density of the standard 
gaussian law $N(0,\mathbf 1_d)$. 
Hence 
\[ 
g(y) = \int_{\mathbb R^d} f(x) \, \frac{1}{\bigl(\sqrt{2 \pi}\,\bigr)^d} 
\, e^{-\|y - x\|^2/2} \, dx \,. 
\] 
Let $\mathcal F$ denote the Fourier transform used by Probabilists: 
\[ 
\mathcal{F}(\varphi)(u) = \int e^{i \langle u, v \rangle} 
\, \varphi(v) \, dv \,, 
\] 
for which 
\[ 
\mathcal{F}(N(0,\mathbf 1_d))(u) = e^{-\|u\|^2/2} 
\] 
and 
\[ 
[\mathcal{F}^{-1}(\psi)](y) = 
\frac{1}{(2 \pi)^d} \int e^{-i\langle y,u\rangle} \, \psi(u)\, du \,. 
\] 
The Fourier transform of a convolution product being the 
product of the transforms, \eqref{blur} implies 
\[ 
\mathcal{F}(g) = \mathcal{F}(f) \, e^{-\|.\|^2/2} \,. 
\] 
Hence 
\begin{equation}\label{inverse} 
f = 
\mathcal{F}^{-1} \bigl(\mathcal{F}(g) \, e^{\|.\|^2/2} \bigr) \,. 
\end{equation} 
Thus \emph{theoretically} one can, using Fourier, 
recover the initial signal. 
But any error in the knowledge of $g$ will make \eqref{inverse} unusable: 
the multiplicative factor $e^{\|.\|^2/2}$ may have considerable effects. 

I doubt that $h$ may have any inverse for the convolution product. 
The right-hand side of \eqref{inverse} has value at $x$ 
\[ 
\frac{1}{(2 \pi)^d} \int e^{-i\langle x,t\rangle} e^{\|t^2\|/2} 
\Bigl[\int e^{i\langle t,y\rangle} \, g(y) \, dy \Bigr] \, dt \,. 
\] 
I do not see how one could eliminate $t$ and get an expression such as 
\[ 
\int k(x-y) \, g(y) \, dy \,. 
\]

\providecommand{\bysame}{\leavevmode\hbox to3em{\hrulefill}\thinspace} 
\providecommand{\MR}{\relax\ifhmode\unskip\space\fi MR } 
\providecommand{\MRhref}[2]{%
 \href{http://www.ams.org/mathscinet-getitem?mr=#1}{#2} 
} 
\providecommand{\href}[2]{#2}

\end{document}